\documentclass[a4paper]{article}

\usepackage[english]{babel}
\usepackage[utf8x]{inputenc}
\usepackage[T1]{fontenc}

\usepackage[top=3cm,bottom=2cm,left=3cm,right=3cm,marginparwidth=1.75cm]{geometry}
\usepackage{authblk}
\usepackage{blindtext}
\usepackage{amsmath, amsthm, amssymb, esint, mathtools, mathabx}
\usepackage{amsfonts}
\usepackage{graphicx}
\graphicspath{{/texfiles/}}
\usepackage[colorinlistoftodos]{todonotes}
\usepackage[colorlinks=true, allcolors=blue]{hyperref}
\usepackage{multicol}

\DeclareMathOperator*{\esssup}{ess\,sup}

\DeclareMathOperator*{\essinf}{ess\,inf}

\numberwithin{equation}{section}

\title{Sparse bounds on variational norms along monomial curves}
\author{A. Martina Neuman}
\affil{Department of Mathematics, New York University, Shanghai}
\affil[ ]{\textit{marsneuman@nyu.edu}}

\begin{document}

\maketitle

\begin{abstract} \noindent Consider a monomial curve $\gamma:\mathbb{R}\to\mathbb{R}^{d}$ and a family of truncated Hilbert transforms along $\gamma$, $\mathcal{H}^{\gamma}$. This paper addresses the possibility of the pointwise sparse domination of the $r$-variation of $\mathcal{H}^{\gamma}$ - namely, whether the following is true:
\begin{equation*}V^{r}\circ\mathcal{H}^{\gamma}f(x)\lesssim \mathcal{S}f(x)\end{equation*}
where $f$ is a nonnegative measurable function, $r>2$ and $\mathcal{S}f(x) = \sum_{Q\in\mathcal{Q}}\langle f\rangle_{Q,p}\chi_{Q}(x)$ for some $p$ and some sparse collection $\mathcal{Q}$ depending on $f,p$.
\end{abstract}


\tableofcontents

\section{Introduction}

\noindent
Consider the following monomial curve
\begin{equation} \label{eq 1.1} \gamma_{P}(t) = (t,t^2,\cdots,t^{d}). \end{equation}
Let $\mathcal{I}$ be a countable set and consider the family $\mathcal{A} = \{A_{t}\}_{t\in\mathcal{I}}$. Let $r>0$. Then the $r$-variation of $\mathcal{A}$ is, 
\begin{equation*}V^{r}\mathcal{A}=V^{r}\{A_{t}\}_{t\in\mathcal{I}} = \sup_{t_{i_1}<\cdots<t_{i_{N}}}\left(\sum_{j=1}^{N-1} |A_{t_{i_{j+1}}} - A_{t_{i_{j}}}|^{r}\right)^{1/r},\end{equation*}
where the supremum runs over all finite increasing subsequences of indices in $\mathcal{I}$. This paper focuses on such $r$-variation of the family of truncated Hilbert transforms along a monomial curve $\gamma$.\\

\noindent
Let $f:\mathbb{R}^{d}\to\mathbb{R}$ be a measurable function. The Hilbert transform of $f$ along a curve $\gamma$ is defined as follows,
\begin{equation*}H^{\gamma}f(x) = pv\int_{-\infty}^{\infty}f(x-\gamma(t))\,\frac{dt}{t}.\end{equation*}
The corresponding family of truncated Hilbert transforms along $\gamma$ is then, $\mathcal{H}^{\gamma} = \{H^{\gamma}_{s}\}_{s>0}$ where $H^{\gamma}_{s}f(x) = \int_{|t|>s}f(x-\gamma(t))\,\frac{dt}{t}$. Define the $r$-variation of the truncated Hilbert transforms of $f$ along $\gamma$, using a set $\mathcal{I}$ (of positive reals), as follows,
\begin{equation} \label{eq 1.2} V^{r}\{H^{\gamma}_{s}f(x)\}_{s\in\mathcal{I}} = \sup_{s_{i_1}<\cdots<s_{i_{N}}}\left(\sum_{j=1}^{N-1} |H^{\gamma}_{s_{i_{j+1}}}f(x) - H^{\gamma}_{s_{i_{j}}}f(x)|^{r}\right)^{1/r},\end{equation}
where the supremum again runs over all finite subsequences $s_{i_1}<\cdots<s_{i_{N}}$. \\

\noindent
For instance, if $\mathcal{I} = \{2^{l}\}_{l\in\mathbb{Z}}$, the sum in $\ref{eq 1.2}$ is simply, $\left(\sum_{j=1}^{N-1} \bigg|\int_{2^{i_{j}}<|t|\leq 2^{i_{j+1}}} f(x-\gamma(t))\,\frac{dt}{t}\bigg|^{r}\right)^{1/r}.$\\

\noindent
Define 
\begin{equation} \label{eq 1.3} \mathcal{T}f(x):=V^{r}\{H^{\gamma}_{s}f(x)\}_{s\in\mathcal{I}}\end{equation}
with $r>2$ and $\gamma$ satisfying $\ref{eq 1.1}$. This paper is devoted to a result on sparse domination in bilinear form of such $\mathcal{T}$. In order to formulate the result, some concepts need to be defined.

\subsection{Definitions}

\subsubsection{The cubes}

\noindent
In this exposition, "a cube" in always means a cube in $\mathbb{R}^{d}$ whose sides are parallel to the coordinate axes. Two families of cubes that reflect well the the geometry of the curve $\gamma$ will be utilized, one is the collection of $\gamma$-cubes, and the other - their cousins - the collection of dyadic $\gamma$-cubes. Their definitions, stated below, are largely borrowed from [1].\\

\noindent
{\bf Definition 1.} [1] A $\gamma$-cube $Q\subset\mathbb{R}^{d}$ is a cube whose side-lengths $\vec{l}(Q) = (l_1,\cdots,l_{d}) = (l^{\alpha_1},\cdots,l^{\alpha_{d}})$ for some $l=l_{Q}>0$. \\

\noindent
{\bf Definition 2.} [1] A dyadic $\gamma$-cube $Q\subset\mathbb{R}^{d}$ is a cube whose side-lengths $\vec{l}(Q) = (2^{\lfloor k\alpha_1\rfloor},\cdots,2^{\lfloor k\alpha_{d}\rfloor})$ for some $k\in\mathbb{Z}$. If $Q$ is one such dyadic $\gamma$-cube, define $l_{Q} := 2^{k}$ for the largest possible $k$.\\

\noindent
The collection of all $\gamma$-cubes is denoted by $\mathcal{Q}^{\gamma}$ and the collection of all dyadic $\gamma$-cubes $\mathcal{D}^{\gamma}$. If $Q_0$ be a cube, then $\mathcal{Q}^{\gamma}(Q_0)$ denotes the collection of all $\gamma$-cubes $Q\subset Q_0$ and $\mathcal{D}^{\gamma}(Q_0)$ the collection of all dyadic $\gamma$-cubes $\tilde{Q}\subset Q_0$. \\

\noindent
The $\gamma$-cubes and dyadic $\gamma$-cubes are essentially equivalent in the following sense. If $Q$ is a dyadic $\gamma$-cube then there exists a $\gamma$-cube $\tilde{Q}$ such that $\tilde{Q}\subset Q$ and $l_{Q}\leq 2l_{\tilde{Q}}$, and vice versa.\\

\noindent
The usefulness of working with dyadic $\gamma$-cubes is that, there exist $3^{d}$ universal shifted dyadic $\gamma$-grids:
\begin{equation*}\left\{2^{\lfloor k\alpha_1\rfloor}\bigg[m_1+\frac{j_1}{3},m_1+1+\frac{j_1}{3}\bigg]\times\cdots\times 2^{\lfloor k\alpha_{d}\rfloor}\bigg[m_{d}+\frac{j_{d}}{3},m_{d}+1+\frac{j_{d}}{3}\bigg]\right\}\end{equation*}
where $k\in\mathbb{Z},\vec{m}\in\mathbb{Z}^{d},\vec{j}\in\{0,1,2\}^{d}$. Now let $\mathcal{D}^{\gamma}$ denote a generic grid of the above type - and temporarily allow a slight abuse of notation - then $\mathcal{D}^{\gamma}$ satisfies the following properties [1]:\\

\noindent
1) $\mathcal{D}^{\gamma} = \bigcup_{k}\mathcal{D}^{\gamma}_{k}$ and each "generation" $\mathcal{D}^{\gamma}_{k}$ partitions $\mathbb{R}^{d}$.\\
2) If $Q_1,Q_2\in\mathcal{D}^{\gamma}$ and $l_{Q_1}\leq l_{Q_2}$, then either $Q_1\subset Q_2$ or $Q_1\cap Q_2=\emptyset$.\\
3) If $Q\in\mathcal{D}^{\gamma}$ and $Q^{(1)}\in\mathcal{D}^{\gamma}$ is the smallest dyadic $\gamma$-cube such that $Q\subsetneq Q^{(1)}$, then $l_{Q^{(1)}}\leq C(d,\gamma)l_{Q}$.\\
4) For every $x\in\mathbb{R}^{d}$, there exists a chain $\{Q_{i}\}\subset\mathcal{D}^{\gamma}$ containing $x$ such that $\lim_{i\to\infty}l_{Q_{i}}=0$.\\

\noindent
It's these properties that enables one to have a Calderon-Zygmund decomposition using these dyadic $\gamma$-cubes.\\

\noindent
{\bf Dilation with cubes.} If $Q\subset\mathbb{R}^{d}$ is a cube with side-lengths $\vec{l}(Q) = (l_1,\cdots,l_{d})$, then $\lambda Q$ is another cube with the same center as $Q$ and side-lengths $\vec{l}(\lambda Q)=(\lambda l_1,\cdots,\lambda l_{d})$. One should be careful that if $Q$ is a $\gamma$-cube or a dyadic $\gamma$-cube then $\lambda Q$ might not be a $\gamma$-cube nor a dyadic $\gamma$-cube. 

\subsubsection{Monotonic functions}

\noindent
Order vectors in $\mathbb{R}^{d}$ in the following manner, $\vec{x}\geq\vec{y}$ if $x_{i}\geq y_{i}$. This is not a complete ordering. Consider functions $f:\mathbb{R}^{d}\to\mathbb{R}$ such that exactly one of the following holds,
\begin{equation*} f(\vec{x})\geq f(\vec{y}) \text{ whenever } \vec{x}\geq\vec{y},\end{equation*} 
\begin{equation*} f(\vec{x})\leq f(\vec{y}) \text{ whenever } \vec{x}\geq\vec{y}.\end{equation*} 
These inequalities are thought to hold almost everywhere. Call such functions monotonic functions. These functions are not integrable but locally integrable. If $f$ is a monotonic function and its weak derivative, $Df$ exists, then exclusively, either 
\begin{equation*} (Df(\vec{x})- Df(\vec{y}))\cdot (\vec{x}-\vec{y})\geq 0\end{equation*} 
happens, or
\begin{equation*} (Df(\vec{x})- Df(\vec{y}))\cdot (\vec{x}-\vec{y})\leq 0.\end{equation*} 
happens, for a.e $\vec{x}\geq\vec{y}$. There is a sub-class of these vector fields that satisfy exclusively either of these two inequalities for a.e $\vec{x},\vec{y}$, called "monotonic vector fields" on $\mathbb{R}^{d}$. See [6] for more on these vector fields. 

\subsubsection{Other definitions}

\noindent
Given a full-body $B\subset\mathbb{R}^{d}$, then $|B|$ denotes its full measure, and given $y\in\mathbb{R}^{d}$, then $|y|$ denotes its vector length. Finally, denote $\langle f\rangle_{Q,p}^{p} = |Q|^{-1}\int_{Q} f(x)^{p}\,dx$ and $\langle f\rangle_{Q} = \langle f\rangle_{Q,1}$.\\

\noindent
A family $\mathcal{S}$ of cubes in $\mathbb{R}^{d}$ is said to be $\eta$-sparse, $0<\eta<1$, if for every $Q\in\mathcal{S}$, there exists $E_{Q}\subset Q$ such that $|E_{Q}|\geq\eta|Q|$ and the sets $\{E_{Q}\}_{Q\in\mathcal{S}}$ are pairwise disjoint. See [11] for instance for this well-known definition.

\subsection{Main theorem}

\noindent
{\bf Theorem 1.} Let $d\geq 2$. Let $\mathcal{T}$ be defined as in $\ref{eq 1.3}$ with $\gamma = \gamma_{P}$. Let $f:\mathbb{R}^{d}\to\mathbb{R}_{\geq 0}$ be a monotonic, measurable function. Fix $Q_0\in\mathcal{D}^{\gamma}$. Let $1<p<\infty$. Then there exist $\kappa=\kappa(d,\gamma)$ and a sparse collection $\mathcal{S}$ of $\gamma$-cubes, such that for a.e $x$,
\begin{equation} \label{eq 1.4}
|\mathcal{T}(f\chi_{\kappa Q_0})(x)|\lesssim \sum_{Q\in\mathcal{S}}\langle f\chi_{\kappa^2 Q_0}\rangle_{\kappa^2 Q}\chi_{Q}(x).\end{equation}

\noindent
{\it Remark:} The implicit constants will be made clear in the body of the proof; they are dependent on natural parameters of the dimensions, the $L^{p}$ exponents, the exponents of the curves and various operator norms associated with $\mathcal{T}$. 

\section{General discussion and acknowledgement}

\noindent
The approach of this paper roughly follows the general principle laid out in [11], where in order for one to obtain a pointwise sparse bound of a sub-linear operator $T$, one needs to obtain a control of its operator norm as well as the operator norm of its associated, so-called tail maximal operator. A more relaxed approach of this principle has just been uploaded [12]. However, this apparatus fails to yield a result in the "variation-sparse-low-dimensional" situation. As a matter of fact, any positive pointwise sparse result for variational norm along a monomial curve should lead to a contradiction.\\ 

\noindent
Related works in this "variation-sparse" direction of consist of the papers [4], [7] and [9]. This list is certainly not complete. 

\subsection{Acknowledgement}

\noindent
In the previous version, there is a mistake in controlling the tail maximal operator, which leads to an incorrect pointwise result, which would hold true for a smaller sub-class of functions, regrettably. The author thanks Professor Lacey for pointing that out.

\section{Proof of Theorem 1}

\subsection{The tail maximal operator}

\noindent
Let $T$ be an operator of functions on $\mathbb{R}^{d}$. The tail maximal operator $\mathcal{M}_{T}$ associated with $T$ is defined by,
\begin{equation*}\mathcal{M}_{T}f(x) = \sup_{Q\ni x}\esssup_{\xi\in Q}|Tf\chi_{\mathbb{R}^{d}\setminus\kappa Q}(\xi)|\end{equation*}
where the supremum is taken over all the cubes containing $x$. A localized version, subject to a cube $Q_0$, of such operator, is,
\begin{equation*}\mathcal{M}_{T,Q_0}f(x) = \sup_{Q\ni x; Q\subset Q_0}\esssup_{\xi\in Q}|Tf\chi_{\kappa Q_0\setminus \kappa Q}(\xi)|.\end{equation*}
The constant $\kappa$ will be chosen beforehand. These definitions are originated in [11]. One needs two divergent points for the definitions of the tail maximal operators used here.\\ 

\noindent
{\bf Definition 3.} Let $C$ be the smallest constant of all $C(d,\gamma)$ satisfying the property (3) in {\bf Definition 2}. Let $\kappa := C+1$. The operators $\mathcal{M}_{T},\mathcal{M}_{T,Q_0}$ are called the $\gamma$-tail maximal operators if the cubes $Q, Q_0$ used are all dyadic $\gamma$-cubes with this defined $\kappa$.\\

\noindent
Regardless of the nature of the cubes used, it's clear that $\mathcal{M}_{\mathcal{T},Q_0}f(x)\leq\mathcal{M}_{\mathcal{T}}f(x).$

\subsection{A preliminary result}

\noindent
{\bf Lemma 4.} Let $1<p<\infty$. For every dyadic $\gamma$-cube $Q_0$ and a.e $x\in Q_0$,
\begin{equation}\label{eq 3.1}
|\mathcal{T}(f\chi_{\kappa Q_0})(x)|\lesssim\|\mathcal{T}:L^{p}\to L^{p}\|f(x)+\mathcal{M}_{\mathcal{T},Q_0}f(x).\end{equation}

\noindent
{\it Proof:} The following route is standard; see [11]. Suppose $x\in Int(Q_0)$ is a point of approximate continuity of $\mathcal{T}(f\chi_{\kappa Q_0})$ and a Lebesgue point of $f$. That means, if $B(x,s)$ denotes a ball of radius $s$ centered at $x$, then for every $\epsilon>0$, the sets
\begin{equation*} E_{s}^{\epsilon}(x) = \{y\in B(x,s):|\mathcal{T}(f\chi_{\kappa Q_0})(y)-\mathcal{T}(f\chi_{\kappa Q_0})(x)|\leq\epsilon\}\end{equation*} 
satisfy $\lim_{s\to 0}\frac{|E_{s}^{\epsilon}(x)|}{|B(x,s)|}=1$. Let $Q(x,s)$ be the smallest dyadic $\gamma$-cube centered at $x$ containing $B(x,s)$ and choose $s$ so small that $Q(x,s)\subset Q_0$. Then, from the definition of $\mathcal{M}_{\mathcal{T},Q_0}$, for a.e $y\in E_{s}^{\epsilon}(x)$, \begin{equation*}|\mathcal{T}(f\chi_{\kappa Q_0})(x)|\leq|\mathcal{T}(f\chi_{\kappa Q(x,s)})(y)|+\mathcal{M}_{\mathcal{T},Q_0}f(x)+\epsilon,\end{equation*}
which implies,
\begin{multline*}|\mathcal{T}(f\chi_{\kappa Q_0})(x)|\leq\essinf_{y\in E_{s}^{\epsilon}(x)}|\mathcal{T}(f\chi_{\kappa Q(x,s)})(y)|+\mathcal{M}_{\mathcal{T},Q_0}f(x)+\epsilon\\ \lesssim_{d}\|\mathcal{T}:L^{p}\to L^{p}\|\frac{1}{|E_{s}^{\epsilon}(x)|^{1/p}}\big(\int_{\kappa Q(x,s)}f^{p}\big)^{1/p}+\mathcal{M}_{\mathcal{T},Q_0}f(x)+\epsilon.\end{multline*}
Now let $s\to 0$ and $\epsilon\to 0$ to obtain $\ref{eq 3.1}$. $\qed$

\subsection{Control of the tail maximal operator}

\noindent
One needs the following concepts of of maximal truncation of the Hilbert transform
\begin{equation*} H^{\gamma}_{*}f(x) = \sup_{0<\epsilon<\rho}\bigg|\int_{\epsilon\leq |t|\leq\rho} f(x-\gamma(t))\,\frac{dt}{t}\bigg|,\end{equation*}
and of a single scale average operator along $\gamma$,
\begin{equation*} A^{\gamma}_{\lambda}f(x) = \int_{\lambda/2 <|t|\leq\lambda} f(x-\gamma(t))\,\frac{dt}{t}.\end{equation*}
\noindent
When working with dyadic scales, one can also write, $A^{\gamma}_{j}f(x) = \int_{2^{j}<|t|\leq 2^{j+1}}f(x-\gamma(t))\,\frac{dt}{t}$.\\

\noindent
Let $f$ be as in the hypotheses of {\bf Theorem 1}. Fix $Q_0\in\mathcal{D}^{\gamma}$. Fix $Q\subset Q_0\cap\mathcal{D}^{\gamma}$ and consider $x,\xi\in Q$. Assume for a moment that $\mathcal{I} = \{2^{l}\}_{l\in\mathbb{Z}}$. For the sake of convenience, the notation for the $\gamma$-tail maximal operator remains unchanged.\\

\noindent
Let $j_{Q}$ be the unique integer such that $l_{Q}=2^{j_{Q}}$. Let $j_{M}$ be the unique integer such that the smallest dyadic $\gamma$-cube containing $\kappa Q_0$ has the length $2^{j_{M}}$; such a cube is contained in $\kappa^2 Q_0$. Given a consecutive finite sequence of integers $j_1<\dots<j_{N}$ that contains $j_{Q}, j_{M}$, one has, for $z\in Q$,
\begin{equation*}\displaystyle \sum_{j=j_1}^{j_{N}-1}\bigg|\int_{2^{j}<|t|\leq 2^{j+1}}f\chi_{\kappa Q_0\setminus\kappa Q} (z-\gamma(t))\,\frac{dt}{t}\bigg|^{r}\leq \sum_{j=j_{Q}}^{j_{M}-1}\bigg|\int_{2^{j}<|t|\leq 2^{j+1}}f(z-\gamma(t))\,\frac{dt}{t}\bigg|^{r},\end{equation*}\\
while a simple $l^1$ domination then leads to,
\begin{equation}\label{eq 3.2}
\left(\sum_{j=j_{Q}}^{j_{M}-1}\bigg|\int_{2^{j}<|t|\leq 2^{j+1}}f (z-\gamma(t))\,\frac{dt}{t}\bigg|^{r}\right)^{1/r} \leq\sum_{j=j_{Q}}^{j_{M}-1}\bigg|\int_{2^{j}<|t|\leq 2^{j+1}}f(z-\gamma(t))\,\frac{dt}{t}\bigg|.\end{equation}
These imply, for $x,\xi\in Q$,
\begin{multline}\label{eq 3.3}
\left(\sum_{j=j_1}^{j_{M}-1}\bigg|\int_{2^{j}<|t|\leq 2^{j+1}}f\chi_{\mathbb{R}^{d}\setminus 3Q} (\xi-\gamma(t))\,\frac{dt}{t}\bigg|^{r}\right)^{1/r}\leq\sum_{j=j_{Q}}^{j_{M}-1}\bigg|\int_{2^{j}<|t|\leq 2^{j+1}}f(\xi-\gamma(t))\,\frac{dt}{t}\bigg| \\ \leq\left\{\sum_{j=j_{Q}}^{j_{M}-1}\bigg|\int_{2^{j}<|t|\leq 2^{j+1}} f(\xi-\gamma(t))\,\frac{dt}{t}\bigg|-\sum_{j=j_{Q}}^{j_{M}-1}\bigg|\int_{2^{j}<|t|\leq 2^{j+1}}f(x-\gamma(t))\,\frac{dt}{t}\bigg|\right\} \\ + \sum_{j=j_{Q}}^{j_{M}-1}\bigg|\int_{2^{j}<|t|\leq 2^{j+1}} f(x-\gamma(t))\,\frac{dt}{t}\bigg|.\end{multline}

\noindent
Let $y=x-\xi$. Then $y=(y_1,\cdots,y_{d})$ with $|y_{i}|\leq l_{Q}^{\alpha_{i}}$ - denote this relation as $|y|\leq l_{Q}$, then the first term on the RHS of $\ref{eq 3.3}$ is, $\sum_{j=j_{Q}}^{j_{M}-1}\big(|\tau_{y}A_{j}^{\gamma}f(x)|-|A_{j}^{\gamma}f(x)|\big)$, and,
\begin{equation} \label{eq 3.4}
\sum_{j=j_{Q}}^{j_{M}-1}\big(|\tau_{y}A_{j}^{\gamma}f(x)|-|A_{j}^{\gamma}f(x)|\big)\leq\sum_{j=j_{Q}}^{j_{M}-1}\sup_{y:|y|\leq l_{Q}}|\tau_{y}A_{j}^{\gamma}f(x)-A_{j}^{\gamma}f(x)|.\end{equation} 

\noindent
Now one has,
\begin{equation*}\bigg|\int_{a<|t|\leq b} f(x-\gamma(t))\,\frac{dt}{t}\bigg| = \bigg|\int_{a}^{b} f(x-\gamma(t))\,\frac{dt}{t} - \int_{a}^{b} f(x-\gamma(-t))\,\frac{dt}{t}\bigg|.\end{equation*}
With $f$ being monotonic and $\gamma=\gamma_{P}$, the term inside the bracket on the RHS of the equality above is either all nonnegative or all nonpositive for all intervals $[a,b]$. That means, for the last term in $\ref{eq 3.3}$,
\begin{equation} \label{eq 3.5} \sum_{j=j_{Q}}^{j_{M}-1}\bigg|\int_{2^{j}<|t|\leq 2^{j+1}} f(x-\gamma(t))\,\frac{dt}{t}\bigg|\leq \bigg|\sum_{j=j_{Q}}^{j_{M}-1}\int_{2^{j}<|t|\leq 2^{j+1}} f\chi_{\kappa^2 Q_0}(x-\gamma(t))\,\frac{dt}{t}\bigg|\leq H^{\gamma}_{*}(f\chi_{\kappa^2 Q_0})(x). \end{equation}

\noindent
From $\ref{eq 3.4}, \ref{eq 3.5}$ and the definition of $\mathcal{M}_{\mathcal{T},Q_0}$, one has for $x\in Q$,
\begin{equation} \label{eq 3.6} \mathcal{M}_{\mathcal{T},Q_0}f(x)\lesssim \sup_{j_1<\cdots<j_{N}}\sum_{j=j_{Q}}^{j_{M}-1}\sup_{y:|y|\leq l_{Q}}\big|\tau_{y}A_{j}^{\gamma}f(x)-A_{j}^{\gamma}f(x)\big|+H^{\gamma}_{*}(f\chi_{\kappa^2 Q_0})(x),\end{equation}
where the supremum runs over all finite increasing consecutive subsequences $j_1<\cdots<j_{N}$ that contains $j_{Q}, j_{M}$.

\subsubsection{The trapezoid}

\noindent
For a dimension $d$, let $\Omega(d)$ be the trapezoid with vertices 
\begin{equation*}(0,0), (1,1), \displaystyle\left(\frac{2}{d+1}, \frac{2(d-1)}{d(d+1)}\right), \left(\frac{d^2-d+2}{d(d+1)},\frac{d-1}{d+1}\right).\end{equation*}
Geometrically, it's a trapezoid that lies on the lower side of the line $y=x$ on the plane. See [2] for the origin of this trapezoid. 

\subsubsection{Control of the tail maximal operator, continued}

\noindent
The following previously known results are needed.\\

\noindent
{\bf Theorem A.} [3] Let $1<p<\infty$. If $f$ is compactly supported, then $\|H^{\gamma}_{*}f\|_{p}\lesssim\|f\|_{p}.$\\

\noindent
To bound the first term in $\ref{eq 3.6}$, one notes that, from van der Corput's lemma and Plancherel's theorem, if $|y_{i}|\leq 1$, 
\begin{equation*}\|A^{\gamma}_1-\tau_{y}A^{\gamma}_1:L^2\to L^2\|\lesssim |y|^{\eta_0}\end{equation*}
for some $\eta_0=\eta_0(\gamma)>0$. One also has from [2] the following result.\\

\noindent
{\bf Theorem B.} [2] If $(1/p,1/q)\in\Omega(d)$, then $A^{\gamma}_1$ is of the restricted weak type $(p,q)$. \\

\noindent
In particular, {\bf Theorem B} says that $A^{\gamma}_1$ is of the restricted weak type $(p,p)$. Interpolation between these two estimates, which is allowed by a result in [14], gives, for $1<p<\infty$,
\begin{equation*}\|A^{\gamma}_1-\tau_{y}A^{\gamma}_1:L^{p}\to L^{p}\|\lesssim |y|^{\eta}\end{equation*}
for some $\eta>0$. By a simple change of variables, this last estimate then implies,
\begin{equation}\label{eq 3.7} \|A^{\gamma}_{\lambda}-\tau_{y}A^{\gamma}_{\lambda}:L^{p}\to L^{p}\|\lesssim\big|\left(\frac{y_1}{\lambda^{\alpha_1}},\cdots,\frac{y_{d}}{\lambda^{\alpha_{d}}}\right)\big|^{\eta},\end{equation}
whenever $|y_{i}|<\lambda^{\alpha_{i}}$. See [1] for a similar approach. Now $\ref{eq 3.7}$ implies the following norm bound for the first term in $\ref{eq 3.6}$,
\begin{equation}\label{eq 3.8} \bigg\|\sum_{j=j_{Q}}^{j_{M}-1}\sup_{y:|y|\leq l_{Q}}(\tau_{y}A_{j}^{\gamma}-A_{j}):L^{p}\to L^{q}\bigg\|\leq \sum_{i\geq 0} |(\frac{1}{2^{i\alpha_1}},\cdots,\frac{1}{2^{i\alpha_{d}}})|^{\eta}\lesssim C(\alpha_1,\eta).\end{equation} 
Then {\bf Theorem A}, $\ref{eq 3.6}, \ref{eq 3.8}$ altogether imply that for a fixed dyadic $\gamma$-cube $Q_0$, then for $1<p<\infty$,
\begin{equation}\label{eq 3.9} \|\mathcal{M}_{\mathcal{T},Q_0}f\|_{p}\lesssim \|f\chi_{\kappa Q_0}\|_{p}+\|f\chi_{\kappa^2 Q_0}\|_{p}.\end{equation}

\subsubsection{Passing to dyadic powers}
 
\noindent
For the case of a general countable set $\mathcal{I}$ of positive reals, suppose $t_{i_1}<\cdots<t_{i_{N}}$ in $\mathcal{I}$ is an increasing subsequence. If $2^{J}\leq t_{i_1}<\cdots<t_{i_{N}}\leq 2^{J+1}$ for some $J$, then argue as in $\ref{eq 3.2}$, one has,
\begin{multline} \label{eq 3.10}
\left(\sum_{j=1}^{N-1}\bigg|\int_{t_{i_{j}}<|t|\leq t_{i_{j+1}}}g (z-\gamma(t))\,\frac{dt}{t}\bigg|^{r}\right)^{1/r} \leq\sum_{j=1}^{N-1}\bigg|\int_{t_{i_{j}}<|t|\leq t_{i_{j+1}}}g(z-\gamma(t))\,\frac{dt}{t}\bigg|\\ \leq\bigg|\int_{2^{J}<|t|\leq 2^{J+1}}g(z-\gamma(t))\,\frac{dt}{t}\bigg|,\end{multline}
for any monotonic $g$. If $[t_{i_1},t_{i_{N}}]$ includes multiple dyadic powers, then one splits the subsequence further into parts that are strictly included in some dyadic intervals and applies the domination in $\ref{eq 3.10}$ again to these parts. More precisely, if $2^{J_1},\cdots, 2^{J_{K}}$ are the dyadic powers interlacing with $t_{i_1},\cdots, t_{i_{N}}$ such that $2^{J_1}\leq t_{i_1}< t_{i_{N}}\leq 2^{J_{K}}$, then
\begin{multline}\label{eq 3.11}
\left(\sum_{j=1}^{N-1}\bigg|\int_{t_{i_{j}}<|t|\leq t_{i_{j+1}}}g (z-\gamma(t))\,\frac{dt}{t}\bigg|^{r}\right)^{1/r} \leq\sum_{j=1}^{N-1}\bigg|\int_{t_{i_{j}}<|t|\leq t_{i_{j+1}}}g(z-\gamma(t))\,\frac{dt}{t}\bigg|\\ \leq\bigg|\int_{2^{J_1}<|t|\leq 2^{J_{K}}}g(z-\gamma(t))\,\frac{dt}{t}\bigg|.\end{multline}
Hence $\ref{eq 3.11}, \ref{eq 3.10}$ imply that $\ref{eq 3.9}$ still holds for the general case of countable set $\mathcal{I}$.

\subsection{Proof of Theorem 1, continued}

\noindent
It is true that
\begin{equation}\label{eq 3.12} \mathcal{T}(f_1+f_2)\lesssim \mathcal{T}f_1+\mathcal{T}f_2.\end{equation} 
Indeed, consider the inequality $(a+b)^{s}\lesssim_{s} a^{s}+b^{s}$ for $a,b\geq 0$ and $s>0$, which is equivalent to $(1+t)^{s}\lesssim_{s} 1+t^{s}$, where $t\geq 1$ and $s>0$, which in turn holds with the constant $2^{s}$, by a simple calculus. Take an increasing subsequence $t_1<\cdots<t_{N}$, then,
\begin{multline*}\left(\sum_{i=1}^{N-1}\bigg|\int_{t_{i}<|t|\leq t_{i+1}} f_1(x-\gamma(t))+f_2(x-\gamma(t))\,\frac{dt}{t}\bigg|^{r}\right)^{1/r} \\ \leq\left(\sum_{i=1}^{N-1}2^{r}\bigg|\int_{t_{i}<|t|\leq t_{i+1}} f_1(x-\gamma(t))\,\frac{dt}{t}\bigg|^{r}+\sum_{i=1}^{N-1}2^{r}\bigg|\int_{t_{i}<|t|\leq t_{i+1}} f_2(x-\gamma(t))\,\frac{dt}{t}\bigg|^{r}\right)^{1/r}\\ \leq 2^{1+1/r}\left\{\left(\sum_{i=1}^{N-1}\bigg|\int_{t_{i}<|t|\leq t_{i+1}} f_1(x-\gamma(t))\,\frac{dt}{t}\bigg|^{r}\right)^{1/r}+\left(\sum_{i=1}^{N-1}\bigg|\int_{t_{i}<|t|\leq t_{i+1}} f_2(x-\gamma(t))\,\frac{dt}{t}\bigg|^{r}\right)^{1/r}\right\}.\end{multline*}\\

\noindent
Assume now that for $1<p<\infty$,
\begin{equation}\label{eq 3.13} \|\mathcal{T}:L^{p}\to L^{p}\|<\infty.\end{equation}
This will be proved later in {\bf Section 4}. For the following {\bf Lemma 5}, let $\mathcal{D}^{\gamma}$ denote one of the finitely many dyadic $\gamma$-grids. The approach is fairly routine.\\

\noindent
{\bf Lemma 5.} There exists a $1/2$-sparse family $\mathcal{F}\subset\mathcal{D}^{\gamma}(Q_0)$ such that for a.e $x\in Q_0$ 
\begin{equation}\label{eq 3.14}
|\mathcal{T}(f\chi_{\kappa Q_0})(x)|\lesssim \sum_{Q\in\mathcal{F}}\langle f\chi_{\kappa^2 Q_0}\rangle_{\kappa^2 Q,p}\chi_{Q}(x).\end{equation}

\noindent
{\it Proof.} Firstly, one wants to show that there exist pairwise disjoint $P_{i}\in\mathcal{D}^{\gamma}(Q_0)$ such that,
\begin{equation}\label{eq 3.15}
|\mathcal{T}(f\chi_{\kappa Q_0})|\chi_{Q_0}\lesssim \langle f\chi_{\kappa^2 Q_0}\rangle_{\kappa^2 Q_0}\chi_{Q_0,p} + \sum_{i}|\mathcal{T}(f\chi_{\kappa P_{i}})|\chi_{P_{i}}.\end{equation}
for a.e $x\in Q_0$ and $\sum_{i}|P_{i}|\leq (1/2)|Q_0|$. To achieve that, let 
\begin{equation*}E=\{x\in Q_0: f(x)>C\langle f\rangle_{\kappa^2 Q_0,p}\}\cup\{x\in Q_0:\mathcal{M}_{\mathcal{T},Q_0}f(x)>C\langle f\rangle_{\kappa^2 Q_0,p}\},\end{equation*} 
with $C$ be sufficiently large so that $|E|\leq\frac{1}{2^{d+2}}|Q_0|$. This is possible due to $\ref{eq 3.9}$. Then one applies the Calderon-Zygmund decomposition to the function $\chi_{E}$ using the dyadic $\gamma$-cubes on $Q_0$ at the height of $\frac{1}{2^{d+1}}$ to produce pairwise disjoint dyadic $\gamma$-cubes $P_{i}\in\mathcal{D}^{\gamma}(Q_0)$ such that 
\begin{center}$|E\setminus\cup_{i}P_{i}|=0 \text{ and } P_{i}\cap E^{c}\not=\emptyset \text{ and } \sum_{i}|P_{i}|\leq (1/2)|Q_0|.$\end{center} 
That means,
\begin{equation}\label{eq 3.16}\forall \text{ a.e } x\in Q_0, f\chi_{Q_0\setminus\cup_{i} P_{i}}\leq C\langle f\rangle_{\kappa^2 Q_0,p}\text{ and }\esssup_{\xi\in P_{i}}|\mathcal{T}(f\chi_{\kappa Q_0\setminus \kappa P_{i}})(\xi)|\leq C\langle f\rangle_{\kappa^2 Q_0,p}.\end{equation}
By $\ref{eq 3.12}$, one has for a.e $x\in Q_0$
\begin{equation}\label{eq 3.17} |\mathcal{T}(f\chi_{\kappa Q_0})|\chi_{Q_0}\leq |\mathcal{T}(f\chi_{\kappa Q_0})|\chi_{Q_0\setminus\cup_{i} P_{i}} + \sum_{i}|\mathcal{T}(f\chi_{\kappa Q_0\setminus \kappa P_{i}})|\chi_{P_{i}} + \sum_{i}|\mathcal{T}(f\chi_{\kappa P_{i}})|\chi_{P_{i}},\end{equation}
of which the first two terms of the RHS are dominated by,
\begin{equation*} |\mathcal{T}(f\chi_{\kappa Q_0})|\chi_{Q_0\setminus\cup_{i} P_{i}} + \sum_{i}|\mathcal{T}(f\chi_{\kappa Q_0\setminus \kappa P_{i}})|\chi_{P_{i}}\lesssim\langle f\rangle_{\kappa^2 Q_0,p},\end{equation*}
through $\ref{eq 3.1}, \ref{eq 3.13}, \ref{eq 3.16}$, which then together with $\ref{eq 3.17}$ leads to $\ref{eq 3.15}$.\\

\noindent
The next steps to get to $\ref{eq 3.14}$ are the iterations of the type of estimate in $\ref{eq 3.15}$ to obtain the sparse family $\mathcal{F} = \{P^{j}_{i}\}_{j\in\mathbb{Z}_{+}}$, where $P^{0} = Q_0$, $\{P^{1}_{i}\}=\{P_{i}\}$ and $\{P^{j}_{i}\}$ are the dyadic $\gamma$-cubes obtained at the $j$th step of the process. $\qed$\\

\noindent
Take a partition of $\mathbb{R}^{d}$ using cubes $R_{l}\in\{\mathcal{D}^{\gamma}_{\vec{j}}:\vec{j}\in\{0,1,2\}^{d}\}$ such that $Q_0\subset\kappa R_{l}$ for each $R_{l}$. Then for each $R_{l}$, generate a sparse family $\mathcal{F}_{l}\subset\mathcal{D}^{\gamma}(R_{l})$ similarly as in $\ref{eq 3.14}$. Then note that if $Q$ is a dyadic $\gamma$-cube with $l_{Q}=2^{k}$ then $3Q\subset\tilde{Q}$ where $\tilde{Q}$ is a $\gamma$-cube concentric with $Q$ whose side-lengths $\vec{l}(\tilde{Q}) = (2^{(k+2)\alpha_1},\cdots,2^{(k+2)\alpha_{d}})$. Let $\mathcal{S}$ be the collection of such $\tilde{Q}$ for each $Q\in\cup_{l}\mathcal{F}_{l}$. Then $\ref{eq 1.4}$ holds for this family $\mathcal{S}$ that is $1/(2\cdot 4^{d-1})$ sparse.

\section{Control of variational norm}

\noindent
The following discussion uses the ideas and notations introduced in [5], [10] and [13].\\

\noindent
Let $1<p<\infty$. It was observed in [10] that in order to obtain $\|\mathcal{T}f\|_{p} = \|V^{r}\{H^{\gamma}_{s}f\}_{s\in\mathcal{I}}\|_{p}\lesssim\|f\|_{p}$, it's sufficient obtain a similar control of its short $2$-variation operator,
\begin{equation} \label{eq 4.1}
\|S_2\{H^{\gamma}_{s}\}_{s\in\mathcal{I}}:L^{p}\to L^{p}\|\lesssim 1.\end{equation}

\noindent
Let $s\in [1,2]$ and $j\in\mathbb{Z}$. Define $\nu_{0,s}, \nu_{j,s}$ respectively by \begin{equation*}\langle\nu_{0,s},f\rangle = \int_{s\leq |u|\leq 2}f(\gamma(u))\,\frac{du}{u},\end{equation*}
and, 
\begin{equation*}\langle\nu_{j,s},f\rangle = \langle\nu_{0,s},f(2^{j}\cdot)\rangle = \int_{s\leq |u|\leq 2}f(\gamma(2^{j}u))\,\frac{du}{u}.\end{equation*}
One observes similarly as in [10] that,
\begin{equation} \label{eq 4.2}
S_2\{H^{\gamma}_{s}f(x)\}_{s\in\mathcal{I}}=(\sum_{j\in\mathbb{Z}}|V_{2,j}(\mathcal{T}f)(x)|^2)^{1/2} = (\sum_{j\in\mathbb{Z}}\|\{\nu_{j,s}\ast f(x)\}_{s\in [1,2]}\|_{v_2}^2)^{1/2}.\end{equation}
From the results in [10], this then entails that a good control of $S_2\{\nu_{j,s}\ast f(x)\}_{s\in [1,2]}$ is sufficient. Let $D=diag(\alpha_1,\cdots,\alpha_{d})$ and $t^{D}x = (t^{\alpha_1}x_1,\cdots, t^{\alpha_{d}}x_{d})$. Denote $E_{j,s}f(x) = \nu_{j,s}\ast f(x)$. Let $\psi$ be a nonnegative smooth function such that $\sum_{k\in\mathbb{Z}}\psi(t2^{-k})\equiv 1$ for all $t>0$. Define the following rescaled truncated Littlewood-Paley operators,
\begin{equation*}\widehat{\Pi_{j,k}f}(\xi) = \psi_{j,k}(\xi)\hat{f}(\xi) := \psi(|2^{jD}2^{-k}\xi|)\hat{f}(\xi).\end{equation*}
Then $E_{j,s}f(x) = \sum_{k\in\mathbb{Z}}\nu_{j,s}\ast\Pi_{j,k}f(x) = \sum_{k\in\mathbb{Z}}f\ast(\nu_{j,s}\ast\check{\psi}_{j,k})(x) :=\sum_{k\in\mathbb{Z}} E_{j,k,s}f(x)$. Denote also that $\tilde{E}_{j,k,s}$  $f(x) = \frac{d}{ds}(\nu_{j,s}\ast\Pi_{j,k}f)(x)$. As discussed - and in view of $\ref{eq 4.2}$ and $l^1$ domination - it's sufficient to obtain the following type of estimate:
\begin{equation} \label{eq 4.3} \|S_2\{E_{j,s}f(x)\}_{s\in [1,2]}:L^{p}\to L^{p}\|\lesssim 1\end{equation}
for $1<p<\infty$ and uniformly for $j\in\mathbb{Z}$. On the other hand, following from [13], one has,
\begin{multline} \label{eq 4.4} S_2\{E_{j,s}f(x)\}_{s\in [1,2]}\lesssim \sum_{k\in\mathbb{Z}} S_2\{E_{j,k,s}f(x)\}_{s\in [1,2]}\\ \lesssim \sum_{k\in\mathbb{Z}} G_{j,k}f(x)^{1/2}\tilde{G}_{j,k}f(x)^{1/2}\lesssim \sum_{k\in\mathbb{Z}}G_{j,k}f(x)+\sum_{k\in\mathbb{Z}}\tilde{G}_{j,k}f(x)\end{multline}
where $G_{j,k}f(x)^2 = \int_1^2 |E_{j,k,s} f(x)|^2\,\frac{ds}{s}$ and $\tilde{G}_{j,k}f(x)^2 = \int_1^2 |\tilde{E}_{j,k,s}f(x)|^2\,\frac{ds}{s}$. This entails,
\begin{multline*} (\sum_{j\in\mathbb{Z}}S_2\{E_{j,s}f(x)\}_{s\in [1,2]}^2)^{1/2}\lesssim (\sum_{j\in\mathbb{Z}}(\sum_{k\in\mathbb{Z}}G_{j,k}f(x))^2)^{1/2}+(\sum_{j\in\mathbb{Z}}(\sum_{k\in\mathbb{Z}}\tilde{G}_{j,k}f(x))^2)^{1/2}\\
= \left(\sum_{j\in\mathbb{Z}}\bigg(\sum_{k\in\mathbb{Z}}\int_1^2 |\nu_{j,s}\ast\Pi_{j,k}f(x)|^2\,\frac{ds}{s}\bigg)\right)^{1/2}+ \left(\sum_{j\in\mathbb{Z}}\bigg(\sum_{k\in\mathbb{Z}}\int_1^2 |\frac{d}{ds}(\nu_{j,s}\ast\Pi_{j,k}f(x)|^2\,\frac{ds}{s}\bigg)\right)^{1/2}.\end{multline*}
Then based on the Littlewood-Paley theory and this last display, it's also sufficient to obtain the estimates of the following type, uniformly for $s\in [1,2]$:
\begin{equation} \label{eq 4.5}
\|\big(\sum_{j\in\mathbb{Z}} |\nu_{j,s}\ast\Pi_{j,k}f|^2\big)^{1/2}\|_{p}\lesssim\|f\|_{p}\end{equation}
and 
\begin{equation} \label{eq 4.6}
\bigg\|\big(\sum_{j\in\mathbb{Z}} |\frac{d}{ds}(\nu_{j,s}\ast\Pi_{j,k}f)|^2\big)^{1/2}\bigg\|_{p}\lesssim\|f\|_{p}.
\end{equation}

\noindent
Following are two approaches in proving $\ref{eq 4.1}$ - the first shows $\ref{eq 4.3}$ and the second shows $\ref{eq 4.5}, \ref{eq 4.6}$.

\subsection{Approach 1}

\noindent
The first approach is analogous to a similar work in [10]. Hence the notations used follow closely to those appeared in [10]. The reader is recommended to refer to the said work for some shortened steps below.\\

\noindent
One first wants to smoothen up the cut-off at $s$ and $2$ for the measures $\nu_{\cdot,s}$'s. Let $\mu$ to be such that $\langle\mu,f\rangle = \int_{1\leq |u|\leq 2} f(\gamma(u))\,\frac{du}{u}$, and consider the smoothened version of $\nu_{0,s}$ as a smoothened, $s^{D}$-dilated, truncated version of $\mu$, as follows:
\begin{equation*} \nu_{0,s} = \nu_{0,s}^{s}+\nu_{0,s}^2=\sum_{m\in\mathbb{Z}_{+}}\nu_{0,s;m}^{s}+\sum_{n\in\mathbb{Z}_{+}}\nu_{0,s;n}^2,\end{equation*} 
where
\begin{equation*}
\begin{split}\langle\nu^{s}_{0,s;m},f\rangle &= \int f(\gamma(su))\theta(|u|)\phi(2^{m}|u-s|)\,\frac{du}{u},\\
\langle\nu^2_{0,s;n},f\rangle &= \int f(\gamma(su))\tilde{\theta}(|u|)\phi(2^{n}|2-u|)\,\frac{du}{u},
\end{split}
\end{equation*}
for some $\theta,\tilde{\theta}$, both smooth functions, and some smooth $\phi$ supported in $(1/2,2)$; all nonnegative. One still has, $\nu_{j,s;m}^{s},\nu_{j,s;n}^2$ as $2^{jD}$-dilated versions of $\nu_{0,s;m}^{s},\nu_{0,s;n}^2$, respectively. Finally, split $E_{j,k,s}$ into $E_{j,k,s;m}^{s}$ and $E_{j,k,s;n}^2$, and $\tilde{E}_{j,k,s}$ into $\tilde{E}_{j,k,s;m}^{s}$ and $\tilde{E}_{j,k,s;n}^2$. Similarly as in [10], observe that,
\begin{equation*}s\frac{d}{ds} \exp(i(\xi_1(2^{j}su)^{\alpha_1}+\cdots+\xi_{d}(2^{j}su)^{\alpha_{d}})) = u\frac{d}{du}\exp(i(\xi_1(2^{j}su)^{\alpha_1}+\cdots+\xi_{d}(2^{j}su)^{\alpha_{d}}))
\end{equation*}
which gives, for $1\leq s\leq u\leq 2$,
\begin{equation*}
\begin{split}
\frac{d}{ds}(\widehat{\nu_{j,s;m}^{s}})((2^{j}s)^{D}\xi) &=\frac{d}{ds}\int \exp(i(\xi_1(2^{j}su)^{\alpha_1}+\cdots+\xi_{d}(2^{j}su)^{\alpha_{d}}))\theta(u)\phi(2^{m}(u-s))\,\frac{du}{u}
\\&= \int\frac{d}{ds}\big\{\exp(i(\xi_1(2^{j}su)^{\alpha_1}+\cdots+\xi_{d}(2^{j}su)^{\alpha_{d}}))\big\}\theta(u)\phi(2^{m}(u-s))\,\frac{du}{u}
\\ &\hspace{30pt}+ \int\exp(i(\xi_1(2^{j}su)^{\alpha_1}+\cdots+\xi_{d}(2^{j}su)^{\alpha_{d}}))\frac{d}{ds}\big\{\theta(u)\phi(2^{m}(u-s))\big\}\,\frac{du}{u}
\\ &= \int\frac{d}{du}\big\{\frac{1}{s}\exp(i(\xi_1(2^{j}su)^{\alpha_1}+\cdots+\xi_{d}(2^{j}su)^{\alpha_{d}}))\big\}\theta(u)\phi(2^{m}(u-s))\,du \\ &\hspace{30pt}+ \int\exp(i(\xi_1(2^{j}su)^{\alpha_1}+\cdots+\xi_{d}(2^{j}su)^{\alpha_{d}}))\frac{d}{ds}\big\{\theta(u)\phi(2^{m}(u-s))\big\}\,\frac{du}{u}
\end{split}
\end{equation*}
which, from integration by parts, in turn gives,
\begin{multline*}
\frac{d}{ds}(\widehat{\nu_{j,s;m}^{s}})((2^{j}s)^{D}\xi) =-\int\frac{1}{s}\exp(i(\xi_1(2^{j}su)^{\alpha_1}+\cdots+\xi_{d}(2^{j}su)^{\alpha_{d}}))\frac{d}{du}\big\{\theta(u)\phi(2^{m}(u-s))\big\}\,du \\
+ \int\exp(i(\xi_1(2^{j}su)^{\alpha_1}+\cdots+\xi_{d}(2^{j}su)^{\alpha_{d}}))\frac{d}{ds}\big\{\theta(u)\phi(2^{m}(u-s))\big\}\,\frac{du}{u}.
\end{multline*}
One has a similar conclusion for the branch $-2\leq u\leq s\leq -1$. A direct application of Van der Corput's lemma yields, for $1\leq s\leq 2$,
\begin{equation} \label{eq 4.7}
\begin{split}
\big|\frac{d}{ds}(\widehat{\nu_{j,s;m}^{s}})((2^{j}s)^{D}\xi)\big|&\lesssim_{\gamma,d} \min\{1,2^{m}|2^{jD}\xi|^{-1/d}\} \\
\big|(\widehat{\nu_{j,s;m}^{s}})((2^{j}s)^{D}\xi)\big|&\lesssim_{\gamma,d} \min\{2^{-m},|2^{jD}\xi|^{-1/d}\}.
\end{split}
\end{equation}
Similarly,
\begin{equation} \label{eq 4.8}
\begin{split}
\big|\frac{d}{ds}(\widehat{\nu_{j,s;n}^2})((2^{j}s)^{D}\xi)\big|&\lesssim \min\{1,2^{n}|2^{jD}\xi|^{-1/d}\} \\
\big|(\widehat{\nu_{j,s;n}^2})((2^{j}s)^{D}\xi)\big|&\lesssim\min\{2^{-n},|2^{jD}\xi|^{-1/d}\}.  
\end{split}
\end{equation}
Furthermore, the cancellation nature of $\nu_{j,s}$'s implies that,
\begin{equation} \label{eq 4.9}
|(\widehat{\nu_{j,s}})((2^{j}s)^{D}\xi)|\lesssim \min\{|2^{jD}\xi|^{-1/d},|2^{jD}\xi|^{1/d}\}\end{equation}
which allows one to obtain for all $1\leq s\leq 2$,
\begin{equation} \label{eq 4.10}
\|(\sum_{j\in\mathbb{Z}}|\nu_{j,s}\ast\Pi_{j,k}f|^2)^{1/2}\|_2\lesssim 2^{-c|k|}\|f\|_2\end{equation}
for some $c=c(\gamma,d)>0$. See also the discussion below in {\bf Subsection 4.3}. Then $\ref{eq 4.7},\ref{eq 4.8}, \ref{eq 4.9}, \ref{eq 4.10}$, the second inequality in $\ref{eq 4.4}$ and Plancherel's theorem imply,
\begin{equation} \label{eq 4.11}
\begin{split}
\|S_2\{E_{j,k,s;m}^{s}f\}_{s\in [1,2]}\|_2 &\lesssim_{\gamma,d} 2^{-cm}2^{-c|k|}2^{-c|j|}\|f\|_2\\
\|S_2\{E_{j,k,s;n}^2f\}_{s\in [1,2]}\|_2 &\lesssim_{\gamma,d} 2^{-cn}2^{-c|k|}2^{-c|j|}\|f\|_2
\end{split}
\end{equation}
for some $c=c(\gamma,d)>0$. \\

\noindent
Denote $\tilde{K}_{j,k,s;m}^{s} = \frac{d}{ds}(\nu_{j,s;m}^{s}\ast\check{\psi}_{j,k})$ and $\tilde{K}_{j,k,s;n}^2 = \frac{d}{ds}(\nu_{j,s;n}^2\ast\check{\psi}_{j,k})$. Then $\tilde{K}_{j,k,s;m}^{s}, \tilde{K}_{j,k,s;n}^2$ are the kernels of $\tilde{E}_{j,k,s;m}^{2},\tilde{E}_{j,k,s;n}^2$ respectively. Let $\rho$ be a metric that is homogeneous with respect to $\{t^{D}\}$ and $C_0>1$ sufficiently large (see [8]). Then it follows from the work in [10] (in particular, the proof of {\it Theorem 1.5}) as well as in [13] that, 
\begin{equation} \label{eq 4.12}
\begin{split}
\int_1^2\int_{\rho(x)\geq C\rho(y)}|\tilde{K}_{j,k,s;m}^{s}(x-y) - \tilde{K}_{j,k,s;m}^{s}(x)|\,dx\frac{ds}{s}&\lesssim |k|,\\
\int_1^2\int_{\rho(x)\geq C\rho(y)}|\tilde{K}_{j,k,s;n}^2(x-y) - \tilde{K}_{j,k,s;n}^2(x)|\,dx\frac{ds}{s}&\lesssim |k|. 
\end{split}
\end{equation}
See also the proof of $\ref{eq 4.6}$ in the {\bf Approach 2} below. Still according to [10], $\ref{eq 4.12}$ then in turn implies, for $1<p<\infty$
\begin{equation} \label{eq 4.13}
\begin{split}
\|S_2\{E_{j,k,s;m}^{s}f\}_{s\in [1,2]}\|_{p} &\lesssim |k|\|f\|_{p}\\
\|S_2\{E_{j,k,s;n}^2f\}_{s\in [1,2]}\|_{p} &\lesssim |k|\|f\|_{p}.
\end{split}
\end{equation}
Interpolation between $\ref{eq 4.11}, \ref{eq 4.13}$ gives, for $1<p<\infty$
\begin{equation*}
\begin{split}
\|S_2\{E_{j,k,s;m}^{s}f\}_{s\in [1,2]}\|_{p} &\lesssim_{\gamma,d,p} 2^{-cm}2^{-c|k|}2^{-c|j|}\|f\|_{p}\\
\|S_2\{E_{j,k,s;n}^2f\}_{s\in [1,2]}\|_{p} &\lesssim_{\gamma,d,p} 2^{-cn}2^{-c|k|}2^{-c|j|}\|f\|_{p}
\end{split}
\end{equation*}
for some $c=c(d,p)>0$. Summing these estimates over $k\in\mathbb{Z}; m,n\in\mathbb{Z}_{+}$, one then obtains the following form of $\ref{eq 4.3}$:
\begin{equation*}
\|S_2\{E_{j,s}f\}_{s\in [1,2]}\|_{p} \lesssim_{\gamma,d,p} 2^{-c|j|}\|f\|_{p}.
\end{equation*}

\subsection{Approach 2}

\noindent
The second approach shows $\ref{eq 4.5},\ref{eq 4.6}$. Fix $s\in [1,2]$ and $k\in\mathbb{Z}$. \\

\noindent
For $\ref{eq 4.6}$, one can use the discussion in the proof of {\it Lemma 6.1} in [10] to view that the kernels $\tilde{\mathcal{K}}= \{\tilde{K}_{j,k,s} = \frac{d}{ds}(\nu_{j,s}\ast\check{\psi}_{j,k})\}_{j\in\mathbb{Z}}$ of the convolution operators $\{f\ast \left(\frac{d}{ds}(\nu_{j,s}\ast\check{\psi}_{j,k})\right)\}_{j}$ as having values in the $l^2$ space. Then an $L^{p}$ bound of $\big(\sum_{j\in\mathbb{Z}} \big|\frac{d}{ds}(\nu_{j,s}\ast\Pi_{j,k}f)\big|^2\big)^{1/2}$ can be obtained from an upper bound $M$ of the following modulus of continuity of $\tilde{K}_{j,k,s}$:
\begin{equation*}\sup_{y\in\mathbb{R}^{d}}\int_{\rho(x)\geq C\rho(y)}\|\tilde{\mathcal{K}}(x-y) - \tilde{\mathcal{K}}(x)\|_{l^2}\, dx \leq M\end{equation*}
where $\rho, C$ have the same meaning as in $\ref{eq 4.12}$. It follows from [13] that one can take $M\lesssim \|\nu_{0,1}\|\lesssim 1.$ This implies, for all $1<p<\infty$.
\begin{equation}\label{eq 4.14}
\bigg\|\big(\sum_{j\in\mathbb{Z}} \big|\frac{d}{ds}(\nu_{j,s}\ast\Pi_{j,k}f)\big|^2\big)^{1/2}\bigg\|_{p}\lesssim \|f\|_{p}
\end{equation}
which is $\ref{eq 4.6}$.\\

\noindent
For $\ref{eq 4.5}$, one can follow the above route in a similar manner, by getting a bound for the modulus of continuity of the right kernel. Then one get a similar estimate as in $\ref{eq 4.14}$ for $\big(\sum_{j\in\mathbb{Z}} |\nu_{j,s}\ast\Pi_{j,k}f|^2\big)^{1/2}$, for $1<p<\infty$. Interpolate this estimate with the one from $\ref{eq 4.10}$, one will obtain the following form of $\ref{eq 4.5}$:
\begin{equation*}
\|\big(\sum_{j\in\mathbb{Z}} |\nu_{j,s}\ast\Pi_{j,k}f|^2\big)^{1/2}\|_{p}\lesssim 2^{-c|k|}\|f\|_{p}\end{equation*}
for some $c=c(\gamma,d,p)>0$. 

\subsection{Discussion of $\ref{eq 4.9}, \ref{eq 4.10}$}

\noindent
It was showed in the proof of {\it Corollary 5.1} in [5] that the cancellation nature of a family of measure $\{\sigma_{j}\}_{j\in\mathbb{Z}}$ ($\hat{\sigma}_{j}(0)=0$) that allows
\begin{equation*}
|\hat{\sigma}_{j}(\xi)|\lesssim\min\{|2^{jD}\xi|^{-1/d}, |2^{jD}\xi|^{1/d}\}.\end{equation*}
There, $\langle\sigma_{j},f\rangle = \int_{2^{j}\leq|t|\leq 2^{j+1}}f(\Gamma(t))\,\frac{dt}{t}$, and the considered $\gamma$ belongs to this $\Gamma$ class of curves. Every step of the derivation of the above fact holds for $\nu_{j,s}$ in place of $\sigma_{j}$ for all $1\leq s\leq 2$. This gives $\ref{eq 4.9}$. One can then construct a partition of unity out of $\psi_{j,k}$ as follows. Let $\{\omega_{k}\}_{k\in\mathbb{Z}}$ be a smooth partition of unity on $\mathbb{R}_{+}$ adapted to the intervals $[2^{-k},2^{-k+1}]$. In particular, one wants $0\leq\omega_{k}\leq 1, \sum_{k}\omega_{k}^2(t)\equiv 1 \text{ if } t>0$ and $supp(\omega_{k})\subset [2^{-k-1},2^{-k+1}]$. Finally, for each $j\in\mathbb{Z}$, let $\{\omega^{j}_{k}\}_{k\in\mathbb{Z}}$ be the $2^{jD}$-dilated version of $\{\omega_{k}\}_{k\in\mathbb{Z}}$. Hence one can think of $(\omega^{j}_{k})^2$ as an appropriately smooth cut-off of $\psi_{j,k}$. Then by Plancherel's theorem,
\begin{equation*}
\sum_{j\in\mathbb{Z}}\||\nu_{j,1}\ast\Pi_{j,k}f|\|_2^2\lesssim \sum_{j\in\mathbb{Z}}\int_{\Delta_{k}^{j}}|\widehat{\nu_{j,1}}|^2(2^{jD}\xi)|\hat{f}|^2(\xi)\,d\xi
\end{equation*}
where $\Delta_{k}^{j} = \{2^{-k-1}\leq |2^{jD}\xi|\leq 2^{-k+1}\}$. If $k<0$ then $|2^{jD}\xi|\geq 2^{-k+1}$. Hence by $\ref{eq 4.9}$,
\begin{equation*}
\sum_{j\in\mathbb{Z}}\||\nu_{j,1}\ast\Pi_{j,k}f|\|_2^2\lesssim \sum_{j\in\mathbb{Z}}\int_{\Delta_{k}^{j}}|\widehat{\nu_{j,1}}|^2(2^{jD}\xi)|\hat{f}|^2(\xi)\,d\xi\lesssim 2^{k/d}\|f\|_2^2.
\end{equation*}
Argue similarly for the cases, $k>1$ and $k=0,1$, one obtains, for all $k\in\mathbb{Z}$,
\begin{equation*}
\sum_{j\in\mathbb{Z}}\||\nu_{j,1}\ast\Pi_{j,k}f|\|_2^2\lesssim 2^{-|k|/d}\|f\|_2.
\end{equation*}
The same is true for all other $\nu_{j,s}$. All this then implies $\ref{eq 4.10}$.



\begin{thebibliography}{11}

\bibitem{CO}
L. Cladek and Y. Ou,
\textit{Sparse domination of Hilbert transforms along curves},
Mathematical Research Letters, 2017.

\bibitem{Ch 1}
M. Christ,
\textit{Convolution, curvature, and combinatorics: A case study}, 
Int. Math. Res. Not. 19 (1988), 1033–1048.

\bibitem{Ch 2}
M. Christ,
\textit{Hilbert transforms along curves - I. Nilpotent groups}, 
Ann. of Math., 122 (1985), 575-596.

\bibitem{dPDU}
F. di Plinio, Y. Do and G. Uraltsev,
\textit{Positive sparse domination of variational Carleson operators},
Annali della Scuola normale superiore di Pisa, Classe di scienze, 18 (2017).

\bibitem{DdF}
J. Duoandikoetxea and J. Rubio de Francia,
\textit{Maximal and singular integral operators via Fourier transform estimates}, 
J.L. Invent. Math., 84 (1986), 541-561.

\bibitem{E}
L. Evans,
\textit{Partial differential equations},
Amer. Math. Soc., 1998.

\bibitem{dFSZ-K}
F. Franca Silva and P. Zorin-Kranich,
\textit{Sparse domination of sharp variational truncations},
math.CA arXiv:1604.05506v3.

\bibitem{G-CdF}
J. Garcia-Cuerva and J. Rubio de Francia,
\textit{Weighted norm inequalities and related topics},
North Holland, 1985.

\bibitem{HLP}
T. Hytonen, M. Lacey and C. Perez
\textit{Sharp weighted bounds for the $q$-variation of singular integrals},
Bull. London Math. Soc., 45 (2013) 529–540.

\bibitem{JSW}
R. Jones, A. Seeger and J. Wright,
\textit{Strong variational and jump inequalities in harmonic analysis}, 
Trans. Amer. Math. Soc., 360 (2008), 6711-6742.

\bibitem{L 1}
A. Lerner
\textit{On pointwise estimates involving sparse operators},
New York J. Math., 22 (2015).

\bibitem{LO}
A. Lerner, S. Ombrosi
\textit{Some remarks on the pointwise sparse domination},
math.CA arXiv:1901.00195v1

\bibitem{dF}
J. Rubio de Francia,
\textit{Maximal functions and Fourier transforms}, 
Duke Math. J., 53 (1986), no. 2, 395-404.

\bibitem{SW 1}
E. Stein and S. Wainger,
\textit{An extension of a theorem of Marcinkiewicz and some of its applications}, J. of Math. and Mechanics, 8 (1959), no. 2, 263-284.

\bibitem{SW 2}
E. Stein and S. Wainger,
\textit{Problems in harmonic analysis related to curvature}, Bull. Amer. Math. Soc., 84 (1978), no. 6, 1239-1295.

\end{thebibliography}
\end{document}